\documentclass[12pt]{article}

\usepackage{amsfonts}
\usepackage{amsmath}
\usepackage{amsthm} %removed this on 12th April because of latex error.
\usepackage{epsfig}
\usepackage{graphics}
\usepackage{pslatex}

\usepackage{latexsym}
\usepackage{amssymb,exscale}
\usepackage[all]{xy}
\usepackage{hyperref}
\usepackage{verbatim} % for commenting out big chunks of text

\theoremstyle{plain} %documentation says there are only three styles
    \newtheorem{theorem}{Theorem}
    \newtheorem{lemma}[theorem]{Lemma}
    
    \newtheorem{corollary}[theorem]{Corollary}

\theoremstyle{definition} % For roman text in the body

    \newtheorem{remark}[theorem]{Remark}

\def\lam{\lambda}

\def\ra{\rightarrow}

\def\C{\mathbb{C}}

\def\R{\mathbb{R}}

\def\intt{\int\limits}

\def\l{\left}
\def\r{\right}
\def\<{\langle}
\def\>{\rangle}

\def\bar{\overline}

\newcommand\tr{{\mbox{\rm tr}}}

\newcommand\mnote[1]{} %off
\newcommand\be{\begin{equation*}}

\newcommand\ee{\end{equation*}}

\newcommand\ben{\begin{equation}}
\newcommand\een{\end{equation}}
\newcommand\bes{\begin{eqnarray*}}
\newcommand\ees{\end{eqnarray*}}

\newcommand{\supp}{\mbox{\rm supp}}

\newcommand{\sm}{{\raise0.3ex\hbox{$\scriptstyle \setminus$}}}

\def\l{\left}
\def\r{\right}

\def\lam{\lambda}

\renewcommand{\phi}{\varphi}

\def\CHI{\mathchoice%
{\raise2pt\hbox{$\chi$}}%
{\raise2pt\hbox{$\chi$}}%
{\raise1.3pt\hbox{$\scriptstyle\chi$}}%
{\raise0.8pt\hbox{$\scriptscriptstyle\chi$}}}
\def\smalloplus{\raise1pt\hbox{$\,\scriptstyle \oplus\;$}}

\author{Alice Guionnet\thanks{
UMPA, CNRS  UMR 5669, ENS Lyon, 46 all\'ee d'Italie,
69007 Lyon, France. {aguionne@umpa.ens-lyon.fr}.
This work was partially supported by the ANR project
ANR-08-BLAN-0311-01.}{ } \;
{ }and Ofer Zeitouni\thanks{School of Mathematics,
University of Minnesota and Faculty of Mathematics,
Weizmann
Institute, POB 26, Rehovot 76100, Israel. {zeitouni@math.umn.edu}. 
The work of this author was partially
supported by NSF grant DMS-0804133 and by a grant from the Israel 
Science Foundation.}}
\title{Support convergence in the single ring theorem}

\begin{document}
 \bibliographystyle{abbrv}
\date{ {\small
December 8, 2010}}
\maketitle

\noindent
{\bf Keywords} Random matrices,   non-commutative measure,
Schwinger--Dyson equation.

\noindent{\bf 
 Mathematics Subject of Classification :} {15A52 (46L50,46L54)}

\begin{abstract}
%{\tt To modify} 
We study  the eigenvalues of  non-normal square matrices
of the form $A_n=U_nT_nV_n$ with $U_n,V_n$ independent Haar distributed 
on the unitary group and $T_n$ real diagonal. We show that when the
 empirical measure of the eigenvalues of $T_n$ converges, and $T_n$
satisfies some technical conditions,  all
 these eigenvalues lie in a single ring.

\end{abstract}

\section{The problem}
In \cite{GKZ}, M. Krishnapur and the authors
considered the convergence of the empricial 
measure of (complex) eigenvalues of matrices of the form $A_n=T_nU_n$,
where $U_n$ is Haar distributed on $\mathcal{U}(n)$, the unitary
group of $n\times n$ matrices, and independent of 
the self-adjoint matrix $T_n$ (which therefore
can be assumed diagonal, with real non-negative entries $s_i^{(n)}$). 
That is, with $\lambda_i^{(n)}$ denoting the eigenvalues of $A_n$,
$L_{A_n}=n^{-1}\sum_{i=1}^n \delta_{\lambda_i^{(n)}}$ 
their empirical measure, and with $L_{T_n}$ the empirical measure 
of the entries of $T_n$,
%and
%$$L_{S_{n}}=\frac1{2n} \sum_{i=1}^n [\delta_{s^{(n)}_i}+
%\delta_{-s^{(n)}_i}]\,$$ 
%the {\it symmetrized} empirical measure of the entries of $T_n$, 
the following is part of the
%We write $G_{T_n}$ for $G_{L_{S_n}}$.
%For a measure $\mu$ supported on $\R_+$, we write 
%
%The 
main result of \cite{GKZ}.
Throughout, for a probability measure $\mu$ supported
on $\R$ or on $\C$, we write
$G_\mu$ for its Stieltjes transform, that is
$$G_\mu(z)=\int \frac{\mu(dx)}{z-x}\,.$$
$G_\mu$ is analytic off the support of $\mu$.
We write $G_{T_n}$ for $G_{\tilde L_{T_n}}$, where for any probability measure
$\mu$ on $\R$
we use $\tilde \mu$ to denote the symmetrized of $\mu$, i.e.
the probability measure satisfying $\tilde \mu(A)=(\mu(A)+\mu(-A))/2$.

\begin{theorem}
	\label{main-theo}
	Assume $\{L_{T_n}\}_n$ converges weakly
	to a probability measure $\Theta$ compactly supported
	on $\R_+$.
	Assume further the following. 
	\begin{enumerate}
		\item
			There exists a constant $M>0$
	so that 
	\begin{equation} \label{bound}\lim_{n\to\infty} P(\|T_n\|>M)=0\,.
	\end{equation}
\item There exist a sequence of events $\{{\cal G}_n\}$ with
	$P({\cal G}_n^c)\to 0$ and
	%OOO0927
	constants $\delta,\delta'>0$ so that for Lebesgue almost any $z\in \C$,
	with $\sigma_n^z$ the minimal singular value
	of $zI-A_n$,
	\begin{equation}
		\label{eq-300609b}
		E({\bf 1}_{{\cal G}_n}
		{\bf 1}_{\{\sigma_n^z<n^{-\delta}\}}
		(\log \sigma_n^z)^2)<\delta'\,.
	\end{equation}
%\item There exist  constants $K,\kappa,\kappa'>0$ such that, for all $n$
%	large,
%	\begin{equation}
%		\label{eq-stboundd}
%	\lim_{n\to\infty}P\left(\sup_{\Im(z)>n^{-\kappa'}}
% 	|G_{\tilde \Theta}(z)-G_{D_n}(z)|\leq 
%		\frac{K}{n^\kappa \Im(z)}
%		\right)=1\,.
%	\end{equation}
\item There exist  constants $\kappa,\kappa_1>0$ such that
	\begin{equation}
		\label{eq-denbound}
		%|\Im G_{\tilde \Theta}(z)|\leq \kappa_1 \quad\mbox{\rm on}\quad \C^+\,.
		|\Im G_{T_n}(z)|\leq \kappa_1 \quad\mbox{\rm on}\quad 
		\{z:\Im(z)>n^{-\kappa}\}\,.
	\end{equation}
	\end{enumerate}
	%OOO0110
	Then 
	%the following hold.
%\begin{itemize}
%\item[a.] 
$L_{A_n}$ converges in probability to a limiting 
	%OOO0930a 
	probability measure
	$\mu_A$,
	 rotationally invariant in $\C$ and supported on the annulus
	$\{re^{i\theta}: a\leq r\leq b\}\,,$
	where  $a=1/\sqrt{\int x^{-2} \Theta(dx)}$ and $b=\sqrt{\int
	x^2\Theta(dx)}.$
\end{theorem}
The conditions of Theorem \ref{main-theo}
	were then showed to hold in some examples of interest, and in
	particular to provide a rigorous proof of the Feinberg-Zee
	``single ring theorem'', see \cite{FZ}. A version of Theorem 
	\ref{main-theo} was also proved to hold when the Haar measure on 
$\mathcal{H}_n$ was replaced by the Haar measure on 	
	the orthogonal group,
	see \cite[Theorem 18]{GKZ}.

Our goal in this paper is to improve the convergence statement in Theorem
\ref{main-theo} to a statement concerning the convergence of
the support of $L_{A_n}$. The following is our main theorem. 
\begin{theorem}
\label{theo-support}
Assume $T_n,U_n$ satisfy the conditions of Theorem \ref{main-theo} and, 
in addition, assume that
\begin{equation}\label{eq-271010a}
a_n:= \frac{1}{\sqrt{\int x^{-2} L_{T_n}(dx)}}\to a=
 \frac{1}{\sqrt{\int x^{-2} \Theta (dx)}}\,,
 \end{equation}
 and
\begin{equation}\label{eq-271010b}
b_n:= {\sqrt{\int x^{2} L_{T_n}(dx)}}\to b=
 {\sqrt{\int x^{2} \Theta (dx)}}\,.
 \end{equation}
 Further assume  if $a>0$ that 
 $\sup_n\|T_n^{-1}\|<\infty$.
 Then, the support of $L_{A_n}$ converges to
 $\supp(\mu_A)=\{z\in \C: |z|\in [a,b]\}$ in probability.
 If moreover the assumptions of Theorem \ref{main-theo} hold
  almost surely with respect to the sequence $T_n$, then
  the convergence
   of the support holds almost surely.
\end{theorem}
When $T_n$ is distributed as the diagonal matrix of singular values
of a Ginibre matrix, the conclusion of Theorem \ref{theo-support} 
follows e.g. from the results in \cite{rider}.

\begin{remark}
	Recall that $\mu_A$ is supported on the annulus $[a,b]\times[0,2\pi)$.
	An elementary computation using the expression for the density
	$\rho_A=\rho_A(r)$ of $\mu_A$,
	 see \cite{GKZ,haageruplarsen}, shows that 
	% its density $\rho_A=
	%\rho_A(r)$ 
	%satisfies
	$$\lim_{r\searrow a} \rho_A(r)=\frac1{\pi a^2} \,,
	%\int x^{-2}\Theta(dx)\,,
	\quad
	\lim_{r\nearrow b} \rho_A(r)=\frac1{\pi b^2}\,.
%	\int x^{2}\Theta(dx)}\,.
	$$
	It is maybe surprising that in spite of the density having a strictly
	positive density at the boundary, the eigenvalues still stick to
	the boundary.
\end{remark}
\subsection{Background and description of the proof}
We recall that
the main difficulty in studying the ESD 
$L_{A_{n}}$ is that $A_n$ is not a normal matrix, that is 
$A_nA^{*}_n\not= A^{*}_nA_n$, almost surely. 
For normal matrices, the limit of 
ESDs can be found by the method of moments or by the method of 
Stieltjes' transforms. For non-normal matrices, the only known method of
proof, which is the one followed
in \cite{GKZ},
is more indirect and follows an idea of Girko~\cite{girko}.
We recall the general outline and some crucial steps which
will be needed in the proof of Theorem \ref{theo-support}.

Introduce the $2n\times 2n$ matrix
\begin{equation}
	\label{eq-hnz}	
 H_{n}^{z}:=\l[\begin{array}{cc} 0 & zI-A_n \\ (zI-A_n)^* & 0 \end{array}  \r].
\end{equation}
Let
$\nu_{n}^{z}$ denote the ESD of $H_{n}^{z}$, 
$$\int \frac{1}{y-x} d\nu_n^z(x)=\frac{1}{2n}\tr\left(
(y-H_{n}^{z})^{-1}
\right)\, ,$$ 
then, see \cite[Eq. (7)]{GKZ},
\begin{equation}\label{eq:girko}
\int \psi(z) dL_{A_{n}}(z) =  \frac{1}{2\pi} \intt_{\C} \Delta \psi(z) 
\int_{\R}\log |x| d\nu_n^z(x)d m(z)\,.\end{equation}
The main advantage of this formulation
is that one can reduce attention to the study 
of the ESD of matrices of the form $(T+U)(T+U)^*$ where $T$ is real diagonal 
 and $U$ is Haar distributed. In the limit (i.e., when $T$ 
and $U$ are replaced by operators in a $C^*$-algebra
that are freely independent,
with $T$ bounded and self adjoint and $U$ unitary), the limit
ESD has been identified by Haagerup and Larsen
\cite{haageruplarsen}. The Schwinger--Dyson equations give both
a characterization of the limit and, more important to us, a discrete
approximation that can be used to estimate the discrepancy between
the pre-limit ESD and its limit. These will play a crucial 
role in the study of the support.

\subsection*{Notation}
We describe our convention concerning constants. Throughout, by the word 
{\em constant} we mean quantities
that are independent of $n$
(or of the complex variables $z$, $z_1$). Generic 
constants
denoted by the letters $C$ or $c$, have
values that may change
from line to line, and they may depend on other parameters. Constants
denoted by $C_i$, $K$, $M$,   $\kappa$ and $\kappa'$
are fixed and do not change from line to line.
\section{Preliminaries:
 evaluation of $\nu^z$ and convergence rates}
\label{sec-aux}
We quickly recall the analysis in \cite{GKZ}, assuming throughout 
that $\|T_n\|$ is uniformly bounded by a constant $M<\infty$.
Fix $z\in \C$ and write $\rho=|z|$.
With  \begin{equation}
     \label{eq-070709a}
   {\bf U}_n=\left(\begin{array}{cc}
 0&U_n\\
 0&0\\
 \end{array}\right)\,,   {\bf Y}_n=\left(\begin{array}{cc}
 0&B_n\\
 B^*_n&0\\
 \end{array}\right),
 \end{equation}
 where $B_n=\rho U_n+T_n$, $T_n$  a real, diagonal
 matrix of uniformly bounded norm
and $U_n$ a ${\mathcal H}_n$ unitary matrix,
define 
$$G^n(z)=E[\frac{1}{2n}\tr\left((z-{\bf Y}_n)^{-1}\right)]\,, \quad
G_{T_n}(z)=G^n(z)|_{\rho=0}$$ 
and
$$G^n_U(z)=E[\frac{1}{2n}\tr\left({\bf U}_n(z-{\bf Y}_n)^{-1}\right)]\,. $$ 
Then, see \cite[Eq. (35)]{GKZ}, the finite $n$ 
Schwinger-Dyson equations for this problem give
\begin{equation}
	\label{eq-120709a}
	\rho (G^n(z_1))^2=2G^n_U(z_1)(1+2\rho G^n_U(z_1))-O_1(n,z_1)\,,
\end{equation}
where 
\begin{eqnarray*} 
O_1(n,z_1)&=&4 E\left[(\frac{1}{2n}\tr- E[\frac{1}{2n}\tr])\otimes 
(\frac{1}{2n}\tr- E[\frac{1}{2n}\tr])\partial
(z_1-{\bf Y}_n)^{-1}{\bf U}_n \right]\\
&=&O\left(\frac{\rho^2}{n^2\Im(z_1)^2 (\Im(z_1)\wedge 1)}\right)\,.
\end{eqnarray*}
In particular,  we have 
\begin{equation}
	\label{eq-070809d}
	G^n_U(z_1)=\frac{1}{4\rho}(-1+\sqrt{1+4\rho^2
	G^n(z_1)^2+4O_1(n,z_1)})\,,
\end{equation}
with  the choice of the square root determined by analyticity 
%OOO011010
and behavior at infinity.
Further, if one defines 
\begin{equation}
	\label{eq-080709c}
	z_2=\psi_n(z_1):=z_1-\frac{\rho^2 G^n(z_1)}{(1+2\rho 
	G^n_U(z_1))}\,,
\end{equation}
then, see \cite[Eq. (39)]{GKZ},
for all $z_1$ with $\Im(z_2)>0$ given by \eqref{eq-080709c},
\begin{equation}\label{cvp0}
G^n(z_1)=
G_{T_n}(\psi_n(z_1))-\tilde O(n,z_1,\psi_n(z_1))\,,
\end{equation}
where
$$\tilde O(n,z_1,z_2)=\frac{2 O(n,z_1,z_2)}{(1+2\rho G_U^n(z_1))}\,$$
and
$$| O(n,z_1,z_2)|\le
\frac{C\rho^2}{n^2 |\Im(z_2)| \Im(z_1)^2 (\Im(z_1)\wedge 1)}\,.$$
In particular, for 
$\Im(z_1)$ 
large,
%(say larger than some $M$), 
it holds that $G^n(z_1)$ and $G_U^n(z_1)$ are small,
implying that $z_2$ is well defined with $\Im(z_2)>0$.
This leads (see \cite[Lemma 10]{GKZ}) to the following weak convergence
statement.

\begin{lemma}
	\label{lem-S-conv}
	If $L_{T_n}$ converges weakly in probability to
	a probability measure $\Theta$, then
	for any
	$z\in \C$,  $\nu_n^z$ converges weakly in probability
	to $\nu^z=\tilde \Theta \boxplus \lam_{|z|}$.
\end{lemma}
\noindent
(Recall that
$\tilde \Theta$
is the symmetrized version of $\Theta$.)

The main work in \cite{GKZ} is then to use the Schwinger-Dyson 
equation \eqref{cvp0} and deduce enough a-priori
bounds that allow one
 to integrate the
logarithmic singularity in \eqref{eq:girko}.
While we will make use of some of these bounds, at this point we
return to our goal, which is to prove Theorem \ref{theo-support}.

\section{Convergence of the support - proof of 
Theorem \ref{theo-support}}
\label{sec-supconv}

Throughout this section, we are in the setup and assumptions of Theorem
\ref{theo-support}.
We first consider the statement concerning convergence
in probability.
Recall that $\supp(\mu_A)=\{z\in \C: |z|\in [a,b]\}$.
Since the density of $\mu_A$ is positive on its support,
see \cite[Remark 8]{GKZ}, we only need to prove that
if $z\not\in \supp(\mu_A)$ then there exists an $\epsilon=\epsilon(z)>0$
so that, with $B(z,\epsilon)$ denoting an open
ball in $\C$ centered at $z$ with radius $\epsilon$,
$$P(L_{A_n}(B(z,\epsilon))\neq 0)\to_{n\to\infty} 0\,.$$

Let $\bar \nu_n^z=\lambda_{|z|}\boxplus \tilde L_{T_n}$
(i.e., $\bar \nu_n^z$ denotes the free convolution of $\lam_{|z|}$
with the symmetrized empirical measure of $T_n$). Since
$L_{T_n}\to \Theta$ weakly, we have that $\bar \nu_n^z\to \nu^z$
weakly. Write $\bar G_n^z$ for the Stieltjes transform
of $\bar \nu_n^z$. Then, $\bar G_n^z(\cdot)$ converges to the Stieltjes 
transform of $\nu^z$, which is denoted in the sequel by $G(\cdot)$.

The first observation we make reduces the study of the support of
$L_{A_n}$ to a question concerning $\bar \nu_n^z$.
\begin{lemma}
\label{lem-sets}
For each $z\not\in \supp(\mu_A)$ there exists an $\epsilon=\epsilon(z)$
so that
$\bar\nu^{z'}_n(B(0,\epsilon))=0$ if $|z-z'|<\epsilon$, for all $n$ large. 
\end{lemma}

Before bringing the proof of Lemma \ref{lem-sets}, we provide 
an a-priori estimate on the spectral radius of certain operators.
Throughout, we use $r(A)$ to denote the spectral radius 
of an operator $A$. We use the convention
that $\|\cdot\|$ denotes 
the operator norm and $\|\cdot\|_2$ the Hilbert-Schmidt norm.
An operator $T$ in a non-commutative probability
space is called {\it $R$-diagonal} 
iff it has the same distribution as $UH$ with $U$ unitary, $H$ positive, and
the algebras generated by $(U,U^*)$ and $H$ freely independent, 
see \cite{haageruplarsen,NS}.
\begin{lemma}
\label{lem-rdiag}
Let $A,B$ be elements of
a non-commutative 
tracial
$C^*$-probability space.
Assume that $A$ is  $R$-diagonal and that there exists a constant
$c_0>0$ so that $\|A\|,\|B\|\leq c_0$. Then, for each $\epsilon>0$ there 
there exists an $\eta=\eta(c_0,\epsilon)>0$ so that
$$ r(A+\eta B)\leq  \|A\|_2+\epsilon\,.$$
\end{lemma}
\noindent
(The case $\eta=0$ of the lemma is 
\cite[Proposition 4.1]{haageruplarsen}.)

\noindent
\proof  
Recall that $r(A+\eta B)=\lim \|(A+\eta B)^n\|^{1/n}$.
By \cite[Corollary 4.2]{haageruplarsen}, we have that $\|A^p\|
\leq (1+p)C \|A\|_2^{p-1}$. Therefore, using the sub-additivity of 
norms, we have, with $C_n=\|(A+\eta B)^n\|$,
\begin{equation}
\label{eq-281010a}
 C_n\leq \|A^n\|+\sum_{k=0}^{n-1}
 \|A^k\|\cdot \|\eta B\|\cdot
 C_{n-k-1}\,,\end{equation}
 where $C_0=1$.

 For $\gamma>0$, set $G(\gamma)=\sum_{n\geq 1} \gamma^n C_n$.
 Clearly $G(\gamma)<\infty$ for $\gamma$ small enough,
 and $r(A+\eta B)^{-1}=\sup \{\gamma: G(\gamma)<\infty\}$.
 Further, $G(\cdot)$ is analytic on $[0,r(A+\eta B)^{-1})$.
 Define also
 $F(\gamma)=
  \sum_{n\geq 1}\gamma^n (1+n) \|A\|_2^{n-1}$ and note that
  $F(\gamma)<\infty$ whenever $\gamma<\|A\|_2^{-1}$.
 From \eqref{eq-281010a} we get
 that whenever $G(\gamma)<\infty$,
 \begin{equation}
 \label{eq-281010b}
 G(\gamma)\leq C \sum_{n\geq 1}\gamma^n (1+n) \|A\|_2^{n-1}
 +|\eta| C c_0
 \sum_{n=1}^\infty \gamma^n \sum_{k=0}^{n-1}
 (1+k)\|A\|_2^{(k-1)\vee 0} C_{n-k-1}\,.
 \end{equation}
 Rearranging, we have that the second sum in the right side of
 \eqref{eq-281010b} equals
 \begin{eqnarray*}
&& \sum_{n=1}^\infty \gamma^n \sum_{k=0}^{n-1}
 (1+k)\|A_2\|^{(k-1)\vee 0} C_{n-k-1}\\
 &=&
 \sum_{k=0}^\infty \|A_2\|^{(k-1)\vee 0} (k+1) \gamma^{k+1}
 \sum_{n=k+1}^\infty
 \gamma^{n-k-1} C_{n-k-1}\\
 &=&
 %\frac{\gamma}{\|A\|_2}
 \gamma\left(1+
 \sum_{k=1}^\infty \|A_2\|^{k-1} (k+1) \gamma^{k}\right)
 (1+G(\gamma))\,.
\end{eqnarray*}
It follows that 
%for some constant $C'=C'(C)$,
$$G(\gamma)\leq 
 CF(\gamma) +C c_0 \eta \gamma (1+F(\gamma)) (G(\gamma)+1)\,.$$
 Therefore, for all $\gamma$ with $G(\gamma)<\infty$ and
 $F(\gamma)<\infty$,
 $$(1-Cc_0\eta \gamma (1+F(\gamma)))G(\gamma)\leq  
 CF(\gamma)+Cc_0\eta\gamma (1+F(\gamma))\,.$$
 It follows that for $\gamma=(\|A\|_2+\epsilon)^{-1}$ 
 there exists an $\eta=\eta(\epsilon,c_0)$ 
 so that  $Cc_0\eta \gamma(1+F(\gamma))<1/2$ and therefore $G(\gamma)<\infty$.
 This implies the statement of Lemma \ref{lem-rdiag}.
 \qed

We can now provide the proof of Lemma 
\ref{lem-sets}.

\noindent
{\it Proof of Lemma \ref{lem-sets}.}
Recall that $\bar\nu^{z'}_n=\tilde L_{T_n}\boxplus \lambda_{|z'|}$, see
Theorem \ref{main-theo}, and thus possesses the same law as 
$X+Y_n$ where $X,Y_n$ are freely independent in a non-commutative
probability space, the law of $X$ is that of  a
Bernoulli $\pm |z'|$ variable, and the law of $Y_n$ being $\tilde L_{T_n}$.
%Recall also that $\mu_A$ is supported on an annulus of inner radius $a=
%1/\|Y^{-1}\|_2$ and outer radius $b=\|Y\|_2$, see Remark 
%\ref{rem-stransf}.

Assume first that $|z|>b$. We may and will assume that for some $\delta>0$,
$|z'|-b_n>\delta>0$ for 
all $n$ large, uniformly in $z'$ with $|z-z'|<\epsilon$,
and consider only such $n$, $\epsilon$ and $\delta$. 
We need to check that there exists an $\epsilon'$ such that
for all $|\eta|<\epsilon'$, $X+Y_n-\eta I$ is invertible.
Writing $X+Y_n-\eta I=X(I+X^{-1}(Y_n-\eta I))$, we see that $X+Y_n-\eta I$ 
is invertible 
iff $I+X^{-1}(Y_n-\eta I)$ is invertible. 
A sufficient condition for that is 
that  $r(X^{-1}(Y_n-\eta I))<1$.
Since $\|X^{-1}\|\leq |z'|^{-1}$ and $\|Y_n\|$ is uniformly bounded,
and since  $X^{-1}Y_n$ is $R$-diagonal
with 
$$\|X^{-1} Y_n\|_2\leq \|X^{-1}\|_2\|Y_n\|_2=
|z'|^{-1} \|Y_n\|_2=|z'|^{-1}b_n\leq 
\zeta<1\,$$
for some fixed $\zeta=\zeta(b,\epsilon,\delta)$, the conclusion follows
from an application of
Lemma \ref{lem-rdiag}   with $A=X^{-1}Y_n$ and $B=X^{-1}$.

Similarly, if $|z|\in [0,a)$ (with $a>0$) and 
%$a>0$ (and thus $Y_n$ is invertible),
$\|Y_n^{-1}\|$ is uniformly bounded,
we repeat the argument, this time writing
$X+Y_n-\eta I=Y_n(I+Y_n^{-1} (X-\eta I))$, and then using
$$\|Y_n^{-1}\|_2\|X\|_2=|z'|/a_n<\zeta<1\,.$$
\qed

Let $${\cal A}=\{z: \exists \epsilon>0, \bar \nu_n^z(B(0,\epsilon))=0
\,,\; \mbox{\rm for all $n$ large}\}.$$
Our next step is to prove a control on $G^n(\cdot)$ for
$z\in {\cal A}$. 
\begin{lemma}
\label{lem-Gconv}
Fix $z\in {\cal A}$, $z\neq 0$. 
Let $\beta>0$ be such that 
for some $n_0$ large enough,
$$[-2\beta,2\beta]\not\in (\cup_{n\ge n_0}
\supp \bar\nu_z^n)\,.$$
%<++>\cup\supp\nu_z\,.$$ 
Then, there are constants $\alpha,\gamma,p>0$ so that
for all $n$ large and 
for all $z_1$ with $\Im(z_1)>n^{-\gamma}$ and 
$\Re(z_1)\in [-\beta,\beta]$,
\begin{equation}
\label{eq-supl2}
|G^n(z_1)-\bar G_n^z(z_1)|<\frac{1}{n^{1+\alpha} \Im(z_1)^p}\,. 
\end{equation}
\end{lemma}
\proof
The proof is divided into several steps. The idea is 
to use (\ref{cvp0}) to compare $G^n$ and $\bar G^z_n$. To do
this up to a small neighborhood of the real
axis, an important point is to show that $G^n$ and $\bar G^z_n$
do not cross the cut of the square root which enters in the definition
of $R_\rho$. The latter point is first shown at
a positive distance of the real axis and then a bootstrap argument is used
to approach the real axis.

\noindent
{\bf Step 1.}
Introduce the set
$${\cal C}_{\epsilon,\beta}= \{z_1:
\Im (z_1)\in [ \epsilon, 2\epsilon), \Re(z_1)\in
[-\beta,\beta]\}.$$
Since $[-\beta,\beta]\not\in \supp \bar \nu^z_n$,
we have that
$\Im(\bar G_n^z(x+i0))=0$ for $x\in [-\beta,\beta]$. 
Moreover $\bar G_n^z$ is uniformly Lipschitz on $
\cup_{\epsilon''\leq \epsilon} {\cal C}_{\epsilon'',\beta}$
(with constant only depending on the distance from $[-\beta,\beta]$
to $\supp \bar\nu_z^n$, which is uniformly bounded below
by $\beta$ by hypothesis).
 Therefore,
for any fixed $\epsilon'(= \beta^{-2} \epsilon)$ (whose value can be 
taken to be $1/12$ in what follows)
we can choose $\epsilon$
small enough such that 
\begin{equation}
\label{eq-patch1}
%\Im(G(z_1))<\epsilon'
\mbox{\rm  for all $z_1\in
\cup_{\epsilon''\leq \epsilon} {\cal C}_{\epsilon'',\beta}$,
it holds that $\Im(\bar G_n^z(z_1))<\epsilon'$, $\Im(G(z_1))<\epsilon'$}.
\end{equation}
By the convergence of $G^n$ to $G$ (which
follows from the weak convergence of $L_{{\bf Y}_n}$ to $\mu_Y$, 
see Lemma \ref{lem-S-conv}),
which can be made uniform by uniform continuity on ${\cal C}_{\epsilon,\beta}$,
and replacing $\epsilon'$ by $3\epsilon'$ if necessary, 
we get  that for all $n>n_0(\epsilon)$,
\begin{equation}
\label{eq-patch2}
\mbox{\rm  for all $z_1\in
{\cal C}_{\epsilon,\beta}$,
it holds that $\Im(G^n(z_1))<3\epsilon'$}.
\end{equation}
%By the convergence of $G^n$ to $G$ (which
%$G^n$ (we use here that $\epsilon$ does not depend on $n$).

\noindent
{\bf Step 2.} Consider $z_1$ with $\Re(z_1)=0$. In that case, the real part 
of both $G^n(z_1)$ and $G(z_1)$ vanishes by symmetry
($G,G^n$ are 
Stieljes transforms of symmetric measures.) 
Now, with $G_U$ as in \cite[Section 3.1]{GKZ},
we have, see \cite[(22)]{GKZ},
$$G_U(z_1)=\frac{1}{4\rho}(-1+\sqrt{1+4\rho^2
	G(z_1)^2})\,.$$
By the analyticity of $ G, G_U$ along the imaginary axis,
we deduce that $\sqrt{1+4\rho^2
	G(z_1)^2}$ can not vanish and since $G(z_1)$ goes to 
zero at infinity,
this implies that $|\Im(G(z_1))|<1/2$. By continuity
for each $\epsilon$ there is a
$\delta=\delta(\epsilon)$ so that
with $z_1$ such that $\Re(z_1)=0, \Im(z_1)>\epsilon$, we have
$|\Im G(z_1)|\le 1/2-\delta$. 
Again by uniform convergence,
%By the convergence of $G^n$ to $G$ (which
%follows from the weak convergence of $L_{{\bf Y}_n}$ to $\mu_Y$, 
%see Lemma \ref{lem-S-conv}),
%which can be made uniform by uniform continuity,
and reducing $\delta$ to $\delta/2$ if necessary, we get the same for
$G^n$ and $\bar G_n^z$.

%\noindent
%{\bf Step 2.}
%Introduce the set
%$${\cal C}_{\epsilon,\beta}= \{z_1:
%\Im (z_1)<\epsilon, \Re(z_1)\in
%[-\beta,\beta]\}.$$
%Since $[-\beta,\beta]\not\in \supp \nu^z$,
%we have that
%$\Im(G(x+i0))=0$ for $x\in [-\beta,\beta]$. Therefore,
%for any fixed $\epsilon'$ (whose value is determined below)
%we can chose $\epsilon$
%small enough such that $\Im(G(z_1))<\epsilon'$ for all $z_1\in
%{\cal C}_{\epsilon,\beta}$.
%Again by uniform convergence,
%this extends to $G^n$.
%

\noindent
{\bf Step 3.}
%Because $G^n$ and $G$ are uniformly bounded on the upper half plane,
%see Lemma \ref{lem-tildegest}, on
Define
$${\cal C}'_{\epsilon,\beta}:=
{\cal C}_{\epsilon,\beta}\cup \{z_1: \Re(z_1)=0, \Im(z_1)>\epsilon\}.$$
By Steps 1 and 2,  there exist $\delta''=\delta''(\epsilon)>0$ such that
\begin{equation}
\label{eq-patch3}
\mbox{\rm  for all $z_1\in
\cup_{\epsilon''\leq \epsilon}
{\cal C}'_{\epsilon'',\beta}$,
it holds that 
%$\Im(G^n(z_1))<2\epsilon'$}.
%for all $n$ large and all
%$z_1\in {\cal C}'_{\epsilon,\beta}$,
$\Re(1+4G^2(z_1))>\delta''$}
\end{equation}
 and, for all $n>n_0(\epsilon)$, 
\begin{equation}
\label{eq-patch4}
\mbox{\rm  for all $z_1\in
{\cal C}'_{\epsilon,\beta}$,
it holds that 
%$\Im(G^n(z_1))<2\epsilon'$}.
%for all $n$ large and all
%$z_1\in {\cal C}'_{\epsilon,\beta}$,
$\Re(1+4(G^n)^2(z_1))>\delta''$}.
\end{equation}
% possess a strictly
%positive real part, which in uniformly bounded away from $0$ there. 
In particular, for all $n>n_0(\epsilon)$,
there is a path leading from $+i\infty$ to any point in
${\cal C}'_{\epsilon,\beta}$
%\cap \{\Im(z_1)>n^{-1/2}\}$, 
along which the choice of the branch of the square-root
in 
%both \eqref{eq3} and 
\eqref{eq-070809d} (and its version with no error term, 
see \cite[Eq. (22)]{GKZ}) is determined by analyticity
(and is the standard one). Denote such a path ${\cal P}$. With this,
we can improve the statement of boundedness in \cite[Lemma 13]{GKZ}
%\ref{lem-tildegest} 
to a convergence statement. In what follows, even though at this
stage the path ${\cal P}$ is bounded away from the real axis
(by $\epsilon$),
we make explicit the dependence of bounds on 
$\Im(z_1)$; this will be useful 
in Step 4.

We rewrite 
\eqref{cvp0} as
\begin{equation}
        \label{eq-071209bnew}
        \tilde G^n(z_1)=G_{T_n}(\psi_n(z_1))=
        G^n(z_1)
        -\tilde O(n,z_1,\psi_n(z_1)).
\end{equation}
With
%$R(z)=(\sqrt{1+4z^2}-1)/2z$ (with the standard choice of cut for 
%the square-root) and 
$$k_n(z_1)=\rho R_\rho(\tilde G^n(z_1))+\psi_n(z_1)-z_1=
\rho R_\rho(\tilde G^n(z_1))-\frac{\rho^2G^n(z_1)}
{(1+2\rho G_U^n(z_1))}\,,$$
%-\frac{G^n(z_1)}{(1+2G^n_U(z_1))},$$
we have 
\begin{equation}
\label{eq-addouf}
 \tilde G^n(z_1)=G_{T_n}\left( z_1+k_n(z_1) -\rho R_\rho (\tilde
G^n(z_1))\right)\,.
\end{equation}
When $\Im(z_1)>0$ is large, we have that $\Im(\psi_n(z_1))$ is large, and
as a consequence, $\tilde G^n(z_1)$ is analytic and small
in this region. It follows that $k_n(z_1)$ is analytic in that
region, and goes to $0$ together with its derivative as 
$\Im(z_1)\to\infty$. Therefore, the map $z_1\to z_1+k_n(z_1)$ is 
invertible in a neighborhood of $+i\infty$ with analytic inverse, denoted
$\varphi_n(z_1)$, which
is a small perturbation of the identity there. Defining
$\hat G^n(z_1)= \tilde G^n(\phi_n(z_1))$, we obtain
$$\hat G^n(z_1)=G_{T_n}(z_1-\rho R_\rho(\hat G^n(z_1)))\,.$$
Comparing with \cite[Equation (29)]{GKZ},
%\eqref{eq6}, 
we get that in a neighborhood of
$+i\infty$, it holds that
$\hat G^n(z_1)=\bar G_n^z(z_1)$, and therefore, in that neighborhood,
%using  \eqref{eq-addouf},
\begin{equation}
\label{eq-addouf1}
\tilde G^n(z_1)=\bar G_n^z(z_1+k_n(z_1)).
\end{equation}
On the other hand,
 from
\eqref{eq-071209bnew},
%and 
%\eqref{eq-stboundd}, 
we have that
\begin{equation}
\label{eq-addoufb}
|\tilde G^n(z_1)-G^n(z_1)|\leq |\tilde O(n,z,\psi_n(z))|
%+\frac{K}{n^\kappa
%\Im(z_1)}
\leq \frac{C\rho^2}{n^2 (\Im(z_1)^4\wedge 1)}\,.
\end{equation}
Thus, for $\Im z \ge C_3 n^{-1/4}$, 
by \eqref{eq-patch4}, \begin{equation}
\label{eq-patch4b}
\mbox{\rm  for all $z_1\in
{\cal C}'_{\epsilon,\beta}$,
it holds that 
$\Re(1+4(\tilde G^n)^2(z_1))>\delta''/2$}.
\end{equation}
Therefore $R_\rho$ is continuously differentiable
at $\tilde G^n(z_1), z_1\in {\cal C}'_{\epsilon,\beta}$
and we have
\begin{equation}\label{ty}
 |\rho R_\rho (\tilde G^n(z_1))-\rho R_\rho(G^n(z_1))|
\leq 
\frac{C}{n^2 (\Im(z_1)^4\wedge 1)}\,.
\end{equation}
Moreover, in the proof of  \cite[Lemma 12]{GKZ},
it was shown that $\rho R_\rho(G^n(z_1))-\frac{\rho^2 G^n(z_1)}{1+2\rho G^n_U(z_1)}$ is small and analytic  on ${\cal C}'_{\epsilon,\beta}$
provided $\epsilon>n^{-1/4}$. Thus, with  \eqref{eq-patch4b},
\eqref{ty},
we deduce that \begin{equation}\label{estkn}
|k_n(z_1)| \leq 
C_{20}/(n^{3/2} (\Im(z_1)^7\wedge 1)
\end{equation}
is smaller than $\Im z_1/2$ and analytic on
${\cal C}'_{\epsilon,\beta}$
provided $\epsilon>n^{-1/7}$.
Hence, \eqref{eq-addouf1} extends to $z_1\in{\cal C}'_{\epsilon,\beta}$
provided $\epsilon>n^{-1/7}$.

Therefore,  again for $z_1\in {\cal C}'_{\epsilon,\beta}$, $\epsilon>n^{-1/7}$,
\begin{eqnarray}
\label{eq-addstep3final}
|G^n(z_1)-\bar G_n^z(z_1)|&\leq & |\tilde G^n(z_1)-\bar 
G_n^z(z_1)|+
|G^n(z_1)-\tilde G^n(z_1)|\nonumber\\
&=& |\bar G_n^z(z_1+k_n(z_1))-\bar G_n^z(z_1)|+|G^n(z_1)-\tilde G^n(z_1)|\nonumber \\
&\leq & \frac{C}{n^{3/2} (\Im(z_1)^8)}\,.
\end{eqnarray}

\noindent
{\bf Step 4} We bootsrap the previous estimate so that one can approach
the real axis: recall that if $S$ denotes the Stieltjes transform
of a probability measure supported on $\R$, we have that
for any $x\in \R$,
$$|\Im(S(x+i\epsilon/2))|\leq 2|\Im(S(x+i\epsilon))|\,.$$
In particular, for all $z_1=x+iy\in {\cal C}_{\epsilon/2,\beta}$, it holds that
\begin{eqnarray*}
|\Im(G^n(z_1))|&\leq &
2|\Im(G^n(x+2iy))|
\\&\leq& 2|\Im(\bar G_n^z(x+2iy))|+ 2|G^n(x+2iy)-\bar G_n^z(x+2iy)|\\
&\leq& 
2\epsilon'+ \frac{2C}{n^{3/2} (\Im(z_1)^8)}\,.
\end{eqnarray*}
In particular, for all $n>n_1(\epsilon)$,
\eqref{eq-patch2} and \eqref{eq-patch4} hold with $\epsilon$ replaced by 
$\epsilon/2$.

One now repeats Step 3, and concludes 
that
\eqref{eq-addstep3final} continues to hold in ${\cal C}'_{\epsilon/2,\beta}$.
Iterating this $\ell$ times
so that $\epsilon
2^{-\ell}\ge n^{-1/7}$(without changing further $n_1(\epsilon)$ or 
$\delta''(\epsilon)$) completes
the proof of Lemma
\ref{lem-Gconv}. \qed

We have the following corollary of
Lemma \ref{lem-Gconv}, whose proof is identical to the 
proof of \cite[Lemma 5.5.5]{AGZ}.
\begin{corollary}
\label{cor-patch1}
With $\beta,\alpha$ as in Lemma \ref{lem-Gconv},
and $\phi$ any smooth function compactly  supported on $[-\beta,\beta]$,
$$\limsup_{n\to\infty} n^{\alpha+1} |E\int \phi d\nu_z^n|<\infty.$$
In particular,
\begin{equation}
	\label{eq-060510a}
	\limsup_{n\to\infty} P(\nu_z^n([-\beta/2,\beta/2])>0)=0\,.
\end{equation}
\end{corollary}
We have now prepared all the steps to prove Theorem \ref{theo-support}.

\noindent
{\bf Proof of Theorem \ref{theo-support}} We only need
to consider $z$ in a compact set. We begin by noting that
\begin{equation}
	\label{eq-050510aa}
P(\mbox{\rm $A_n$ has an eigenvalue in $B(z,\epsilon)$})=
P(\mbox{\rm $\nu_n^{z'}(\{0\})\ge\frac{1}{n}$ for some $z'\in 
B(z,\epsilon)$})\,.
\end{equation}
We write ${\bf Y}_n(z)$ to emphasize the dependence of ${\bf Y}_n$
in $z$. Let
$$\lambda^*({\bf Y}_n(z))= \min\{|\lambda_i({\bf Y}(z))|\}.$$ 
Since ${\bf Y_n}(z)-{\bf Y_n}(z')$ is Hermitian and of norm bounded by
$|z-z'|$,  we 
have that $|\lambda^*({\bf Y}_n(z))-\lambda^*({\bf Y}_n(z'))|\leq |z-z'|$.
%where $f(z,z')$ is continuous (and vanishes on the diagonal).
Thus, for each $z\not\in \supp (\mu_A)$, and with $\beta=\beta(z)$ 
as in Lemma \ref{lem-Gconv}, we can find an $\epsilon=\epsilon(z)$
so that by Chebyshev's inequality
$$P(\nu_n^{z'}(\{0\})\ge\frac{1}{n} \mbox{ for some }z'\in 
B(z,\epsilon))
\leq 
P(\nu_n^{z}([-\beta/2,\beta/2])\ge\frac{1}{n})
\leq Cn^{-\alpha}\to_{n\to\infty} 0\,.$$
Combined with \eqref{eq-050510aa}, we conclude that
$$P(\mbox{\rm $A_n$ has an eigenvalue in $B(z,\epsilon)$})\to_{n\to\infty} 0\,.$$
By a standard covering argument,
this implies that for any compact $G$ with $G\cap (\supp \mu_A)=\emptyset$, it 
holds that
$$P(\mbox{\rm $A_n$ has an eigenvalue in $G$})\to_{n\to\infty} 0\,.$$
This completes the convergence in probability in the statement of
Theorem \ref{theo-support}.

We finally prove the almost sure convergence
by generalizing the ideas of \cite{HT} based
on Poincar\'e inequality. In our case, we shall use concentration of measures
on SU(N) \cite[Theorem 4.4.27]{AGZ}. 
Since we now assume that the assumptions of Theorem \ref{main-theo}
hold for almost all sequence $T_n$, we may and will 
assume the sequence $T_n$ deterministic in the sequel.
Recall that for any bounded measurable function
$\phi$, $\int \phi(x)d\nu_n^z(x)$ is a bounded measurable function
of the random matrix $W_n=U_n^*V_n^*$.
We denote by $E_{ U(n)}$ (resp. $E_{SU(n)}$)
the expectation over $W_n$ following the Haar measure 
on ${\mathcal U}(n)$ (resp. $SU(n)$).
We also write in the sequel ${\cal B}=(\supp\, \mu_A)^c$.

\begin{lemma} 
Fix  $z\in\mathcal B $, $\alpha$ and $\beta$ as in Lemma \ref{lem-Gconv},
and 
a bounded non negative  smooth function $\phi$ with support in $[-\beta,\beta]$.
%vanishing
%which vanishes
%on the support of $\nu^z$.
\begin{enumerate}
\item There exists a finite constant $C$ such
that 
\begin{equation}\label{jkl0}|E_{U(n)}
 [\int \phi(x)d\nu^z_n(x)]|\le \frac{C}{n^{1+\alpha}}\,.\end{equation}
\item  For all $\delta>0$, there exists $z'\in \mathcal B$
so that $|z-z'|\le \delta$ and 
$$|E_{SU(n)}
 [\int \phi(x)d\nu^{z'}_n(x)]|\le  \frac{C}{n^{1+\frac{\alpha}{2}}}\,.$$
Moreover there exists $n_0=n_0(z',\omega)$ so that
for almost every $\omega$ and all $n>n_0$,
%for $n$ large enough, we have
\begin{equation}\label{jkl}
 |\int \phi(x)d\nu^{z'}_n(x)|\le  \frac{1}{n^{1+\frac{\alpha}{16}}}\,.
\end{equation}
\end{enumerate}
\end{lemma}
The last point proves the theorem as  $ A_n$ has
an eigenvalue in $B(z,\epsilon)\subset \mathcal B$
for $\epsilon$ small enough only if 
$$\nu_n^{z'}([-2{\epsilon},2\epsilon])
\ge \frac{1}{n}$$
for all $z'\in B(z,c \epsilon)$, for an appropriate
$c=c(M,z)$.  \eqref{jkl} shows that this is impossible 
 for $n$ sufficiently large,
almost surely. 

\noindent
{\bf Proof.}
The first point of the lemma is a restatement of
the first part of Corollary \ref{cor-patch1}.
%direct consequence 
%of \cite[Lemma 5.5.5]{AGZ} and \eqref{eq-supl2}.
%We next fix $\phi$
For the second, recall that any matrix $W_n$ in the unitary 
group can be decomposed as $W_n= e^{i\theta} S_n$
with $S_n$ in the special unitary group $SU(n)$
and note that multiplying $S_n$ by $e^{i\theta}$ amounts
to rotating $z$ by $e^{i\theta}$ in $H_n^z$.
Therefore, by the Chebyshev
 inequality we deduce from the
first point that
the set $R_n$ of $\theta\in [0,2\pi]$ such that
\begin{equation}\label{bound2}
 |E_{SU(n)}
 [\int \phi(x)d\nu^{e^{i\theta} z}_n(x)]|\le n^{-1-\frac{\alpha}{2}}
\end{equation}
satisfies $|R_n|/2\pi\geq
%has probability  greater than $
1-C n^{-\alpha/2}$, where $|R_n|$ denotes the Lebesgue measure of
$R_n$. Thus, in any interval
of width $n^{-\alpha/2}$ in the circle of radius $|z|$
there is at least an element of $R_n$. We finally
cover the compact set $\mathcal B\cap [0,M]$ 
(with $M$ as in \eqref{bound}) with a covering with mesh $\delta/2$
to obtain the existence of a 
family $(z_i)_{i\ge 0}$
of points of $\mathcal B$ so that 
\eqref{bound2} hold. Repeating this argument
%We may, up to reproduce the
%same argument for 
with the function $ \phi'(x)^2$, 
we also have that
%assume as well that 
\begin{equation}\label{bound3}
 |E_{SU(n)}
 [\int \phi'(x)^2 d\nu^{z_i}_n(x)]|\le C 
 n^{-1-\frac{\alpha}{2}}\,.\end{equation}
Next, remark that $U_n\ra \int \phi(x)d\nu^{z_i}_n(x)$ 
is Lipschitz with constant bounded above by $C\left( n^{-1}
\int \phi'(x)^2 d\nu^{z_i}_n(x)\right)^{\frac 1 2}$.
%, which  is bounded by 
%$n^{-\frac{1}{2} -\frac{\alpha}{4}}$
%with probability greater or equal to $1-n^{-1-\frac{\alpha}{4}}$
%by \eqref{bound3}. 
Set 
$C_n=\{W_n\in SU(n):\int \phi'(x)^2 d\nu^{z_i}_n(x)\le n^{-\frac{\alpha}{4}}\}$.
Then, 
\begin{equation}
	\label{eq-cn}
	P(C_n^c)\leq C n^{-1-\alpha/4}\,.
\end{equation}
Consequently, using \eqref{bound3},
$$E_{SU(n)} [1_{C_n}
\int \phi(x)d\nu^{z_i}_n(x)]\leq C n^{-1-\alpha/4}\,.
$$
Therefore,
we get that for all $n$ large enough,
\begin{eqnarray*}
&&P\left( \left|\int \phi(x)d\nu^{z_i}_n(x)\right|
\ge n^{-1-\frac{\alpha}{16}}\right)\nonumber\\
&\le& P\left( \left|\int \phi(x)d\nu^{z_i}_n(x)-
E_{SU(n)} [1_{C_n}
 \int \phi(x)d\nu^{z_i}_n(x)]
\right|
\ge \frac{1}{2} n^{-1-\frac{\alpha}{16}}\right)\\
&\le&  n^{-1-\frac{\alpha}{4}}\\
&&+
 P\left( \left\lbrace \left|\int \phi(x)d\nu^{z_i}_n(x)-E_{SU(n)} [1_{C_n}
 \int \phi(x)d\nu^{z_i}_n(x)]
\right|
\ge \frac{1}{2} n^{-1-\frac{\alpha}{16}}\right\rbrace\cap C_n\right)\\
&\le& Cn^{-1-\frac{\alpha}{4}}+Ce^{-n^{-2-\frac{\alpha}{8}} n^2 
n^{\frac{\alpha}{2}}}\,,
\end{eqnarray*}
where we have applied  \cite[Theorem 4.4.27]{AGZ} 
to the extension
of the function  $W_n\to g(W_n)=\int \phi(x)d\nu^{z_i}_n(x)$
outside $C_n$ 
which is globally  Lipschitz with constant $n^{-\frac{1}{2}-\frac{\alpha}{4}}$ 
 and uniformly bounded, see e.g. \cite[Section 5.4]{Guionnet} for the existence
 of such extension.
Applying the Borel-Cantelli lemma completes the proof.
\qed

%{\bf Acknowledgments:} We thank Greg Anderson for many
%fruitful and encouraging discussions.
% We thank Yan Fyodorov 
% for pointing out the paper \cite{horn} and
%  Philippe Biane
% for suggesting that our technique could be 
% applied to the examples in \cite{bl01}.
% We thank the referee for a careful reading of the manuscript.
%% respectively.

\noindent
{\bf Acknowledgement} This paper was written while both authors participated 
in the 2010
MSRI program on 
random matrix theory, interacting particle systems and integrable systems.
We thank MSRI for its hospitality.

\bibliographystyle{amsplain}

\end{document}